\documentclass{article}

   \usepackage[dvips]{graphicx}

\usepackage{url}
\usepackage{cite}

\usepackage{graphicx}
\usepackage{amsmath,amsfonts,amssymb}
\usepackage{bm} % bold greek letters
\usepackage{cases}
\usepackage{float}
\usepackage{mwe}
\usepackage{mathtools}
\usepackage{comment}
\usepackage{empheq}
\usepackage{booktabs}
\usepackage{textcomp}
\usepackage{algorithm}%,algcompatible
\usepackage{algorithmic}
\usepackage{array}
\usepackage{subcaption}
\usepackage{graphicx}
\usepackage{url}
\usepackage{caption}
\usepackage{multirow}
\usepackage{systeme}
\usepackage{adjustbox}
\usepackage{soul}
\usepackage{xcolor}
\usepackage{enumitem}
\usepackage[scale=0.8]{geometry}
\usepackage[roman]{written}
\newtheorem{rem}{Remark}

\usepackage{tikz}
\usetikzlibrary{positioning}
\usetikzlibrary{3d}
\usetikzlibrary{matrix}
\usetikzlibrary{decorations.text}
\usetikzlibrary{spy}

\makeatletter
\tikzoption{canvas is plane}[]{\@setOxy#1}
\def\@setOxy O(#1,#2,#3)x(#4,#5,#6)y(#7,#8,#9)%
{\def\tikz@plane@origin{\pgfpointxyz{#1}{#2}{#3}}%
	\def\tikz@plane@x{\pgfpointxyz{#4+#1}{#5+#2}{#6+#3}}%
	\def\tikz@plane@y{\pgfpointxyz{#7+#1}{#8+#2}{#9+#3}}%
	\tikz@canvas@is@plane
}
\makeatother
\tikzoption{canvas is xy plane at z}[]{%
	\def\tikz@plane@origin{\pgfpointxyz{0}{0}{#1}}%
	\def\tikz@plane@x{\pgfpointxyz{1}{0}{#1}}%
	\def\tikz@plane@y{\pgfpointxyz{0}{1}{#1}}%
	\tikz@canvas@is@plane
}
\makeatother

\usepackage[english]{babel}
\usepackage[utf8]{inputenc} % recommended for inputting accented characters
\usepackage[colorlinks=true,linkcolor=black,anchorcolor=black,citecolor=black,filecolor=black,menucolor=black,runcolor=black,urlcolor=black]{hyperref}       % hyperlinks
\usepackage{tipa}
 \expandafter\let\csname equation*\endcsname\relax
 \expandafter\let\csname endequation*\endcsname\relax
 
\usepackage{array}
\makeatletter
\newcommand{\thickhline}{%
    \noalign {\ifnum 0=`}\fi \hrule height 1pt
    \futurelet \reserved@a \@xhline
}
\newcolumntype{"}{@{\hskip\tabcolsep\vrule width 1pt\hskip\tabcolsep}}
\makeatother

\definecolor{red_our}{rgb}{1,0,0}
\definecolor{lime_our}{rgb}{0,1,0}
\definecolor{blue_our}{rgb}{0,0,1}
\definecolor{yellow_our}{rgb}{1,1,0}
\definecolor{cyan_our}{rgb}{0,1,1}
\definecolor{magenta_our}{rgb}{1,0,1}
\definecolor{green_our}{rgb}{0,0.5,0}
\definecolor{orange_our}{rgb}{1,0.6470,0}
\definecolor{acqua_our}{rgb}{0,0.5,1}
\definecolor{blueviolet_our}{rgb}{0.5411,0.1686,0.8862}
\definecolor{deeppink_our}{rgb}{1,0.0784,0.5764}
\definecolor{lemon_our}{rgb}{1,0.9803,0.8039}
\definecolor{midnightblue_our}{rgb}{0.9803,9803,0.4392}
\definecolor{acquamarine_our}{rgb}{0.4980,1,0.9313}
\definecolor{cornflowerblue_our}{rgb}{0.3921,0.5843,0.9294}
\definecolor{maron_our}{rgb}{0.5019,0,0}

\definecolor{green_mat}{rgb}{0.01,0.75,0.24}
\definecolor{purple_mat}{rgb}{0.49,0.18,0.56}
\definecolor{blue_mat}{rgb}{0.00,0.45,0.74}

\begin{document}

\title{Constrained Regularization by Denoising with Automatic Parameter Selection}

\author{Pasquale~Cascarano\thanks{This material is based upon work partially supported by the NSF CAREER award under grant CCF-2043134. This work has been conducted under the activities of INdAM—Gruppo Nazionale per il Calcolo Scientifico. This study was carried out within the MICS (Made in Italy – Circular and Sustainable) Extended Partnership and received funding from the European Union Next-GenerationEU (Piano Nazionale di ripresa e resilienza (PNRR) – Missione 4 Componente 2, Investimento 1.3 - D.D. 1551.11-10-2022, PE00000004). Pasquale Cascarano is with the Department of the Arts, University of Bologna, Bologna, 40123 Italy (e-mail: pasquale.cascarano2@unibo.it).}, Alessandro~Benfenati\thanks{Alessandro Benfenati is with the Department of Environmental Science and Policy, University of Milan, Milan, 20133 Italy (e-mail: alessandro.benfenati@unimi.it).}, Ulugbek~S.~Kamilov\thanks{Ulugbek S. Kamilov is with the Department of Computer Science \& Engineering and Department of Electrical \& System Engineering, Washington University in St. Louis, St. Louis, MO 63130 USA (e-mail: kamilov@wustl.edu).},~\emph{Senior Member,~IEEE}, Xiaojian~Xu\thanks{Xiaojian Xu is with the Department of Electrical Engineering and Computer Science, University of Michigan, Ann Arbor, MI 48109 USA (e-mail: xjxu@umich.edu).}
}

\markboth{}{}
\maketitle

\begin{abstract}
Regularization by Denoising (RED) is a well-known method for solving image restoration problems by using learned image denoisers as priors. Since the regularization parameter in the traditional RED does not have any physical interpretation, it does not provide an approach for automatic parameter selection. This letter addresses this issue by introducing the Constrained Regularization by Denoising (CRED) method that reformulates RED as a constrained optimization problem where the regularization parameter corresponds directly to the amount of noise in the measurements. The solution to the constrained problem is solved by designing an efficient method based on alternating direction method of multipliers (ADMM). Our experiments show that CRED outperforms the competing methods in terms of stability and robustness, while also achieving competitive performances in terms of image quality. %along with residual noise estimators. 
\end{abstract}

\providecommand{\keywords}[1]{\textbf{{Keywords: }} #1}

\section{Introduction}

The problem of recovering an image $\vx \in \mathbb{R}^{n}$ from its degraded measurement $\vb \in \mathbb{R}^{n}$ can be cast as the following linear inverse problem
\begin{equation}\label{inverse-problem}
    \text{find} \quad \vx \in \mathbb{R}^{n} \quad \text{such that} \quad  \vb = \vA \vx + \veta,
\end{equation}
where $\vA \in \mathbb{R}^{n \times n}$ is a known measurement operator and $\veta \in \mathbb{R}^{n}$ is random noise with standard deviation $\sigma_{\veta}$. 

Linear inverse problems are at the core of many applications \cite{cascarano2020super,cascarano2022deepcel0,teodoro2018convergent,mylonopoulos2022constrained,Benfenati_2015}. However, since most inverse problems are ill-posed, it is common to formulate the solution $\vx^{*}\in\mathbb{R}^{n}$ of \eqref{inverse-problem} as a minimizer of a regularized objective function 
\begin{equation}  \label{eq:reg_mod}
    \vx^{*} =  \underset{\vx \in \mathbb{R}^{n}}{\amin} \ \ell (\vx;\vb) + \lambda \rho(\vx).
\end{equation}
The data-fidelity term $\ell(\vx;\vb)$ encodes information on the noise statistics, e.g.\ additive white Gaussian noise (AWGN) assumptions entail $\ell(\vx;\vb)=\frac{1}{2}\lVert \vA \vx - \vb \rVert^{2}_{2}$. 
The regularization parameter $\lambda >0$ is often hand-tuned to obtain optimal restored images. Alternatively, it can be estimated using well established methods such as the discrepancy principle, L-curve, or cross-validation \cite{bertero2006regularization}.
However, the primary challenge consists in designing the regularization functional $\rho(\vx)$ in order to capture the intricate image features. 

Plug-and-Play (PnP) Priors framework has recently emerged as a powerful tool for exploiting sophisticated denoisers as regularizers without explicitly defining $\rho(\vx)$~\cite{venkatakrishnan2013plug,cascarano2022plug,sun2021scalable,kamilov2023plug,kamilov2017plug,sun2019online,xu2020provable}. 
However, the lack of an explicit objective function complicates the analysis of PnP methods in terms of theoretical understanding and convergence guarantees. Regularization by denoising (RED)~\cite{Romano17} is a variant of PnP based on formulating an explicit regularization functional
\begin{equation}
\rho_{\textsf{\tiny RED}}(\vx) := \dfrac{1}{2}\left(\vx^{T} \left(\vx - f(\vx) \right)\right),
\end{equation}
where $f(\cdot)$ denotes a denoiser. When the denoiser is differentiable, locally homogeneous, and has a symmetric Jacobian~\cite{Romano17,reehorst2018regularization}, $\rho_{\textsf{\tiny RED}}$ is convex and its gradient can be efficiently computed as  $\nabla \rho_{\textsf{\tiny RED}}(\vx) = \vx - f(\vx) $. The interpretation with an explicit regularizer simplifies the theoretical analysis of RED algorithms as methods for computing global minimizers of convex objective functions. These conditions are often not satisfied for many practical denoisers~\cite{reehorst2018regularization}, however, the RED algorithms achieve state-of-the-art performances in many imaging applications. 

RED was recently reformulated as a constrained optimization problem by projecting  the least square minimum onto the fixed-point sets of demicontractive denoisers, which are proven to be convex sets~\cite{Elad21}. Ideally, an image denoiser $f$ should satisfy the condition $\mathcal{M} \subset {\rm Fix}(f)$, where $\mathcal{M}$ denotes the manifold of natural images.
In practice, denoisers are often far from being ideal, and their fixed-point sets may not correspond to the set of natural images, leading to suboptimal recovery results.
Additionally, even if any non-expansive denoiser is demicontractive, determining whether a given denoiser is non-expansive or demicontractive is a challenge.

In this letter, we present a constrained RED (CRED) approach inspired by the discrepancy principle \cite{engl1987discrepancy,zanni2015numerical} solved via the alternating direction method of multipliers (ADMM). 
In order to overcome the limits regarding ${\rm Fix}(f)$, we reverse the RED-PRO \cite{Elad21} formulation by considering the RED regularization functional subject to the discrepancy between the measured data and the reconstruction being below a given threshold. This threshold represents the strength of the regularization and has a precise physical significance since it reflects the standard deviation of the noise affecting the data. Different from the regularization parameter of the original unconstrained RED formulation \cite{Romano17}, which must be hand-tuned, the threshold can be estimated by using well-established techniques \cite{immerkaer1996fast}. Therefore, our approach avoids any parameter tuning which may be limiting in real applications. 
Our formulation and the corresponding ADMM scheme are presented in Section \ref{sec:method}. In Section \ref{sec:numexp}, we underline the quality of CRED in terms of image quality metrics and robusteness, with respect to the choice of the denoiser and the ADMM parameters, through several comparisons with RED and RED-PRO.  
\section{Constrained RED}
\label{sec:method}
Upon AWGN assumptions, CRED seeks to compute a solution of \eqref{inverse-problem} by solving the constrained problem
\begin{equation}
    \label{eq:constrained}
        \argmin{\vx \in\mathbb{R}^n}{\rho_{\textsf{\tiny RED}}(\vx)}   \quad\text{subject to}\quad  \lVert \vA \vx - \vb \rVert^{2}_{2} \leq \delta, 
\end{equation}
where $\delta:=\tau \, \sqrt{n} \, \sigma_{\veta}$ with $\tau \in [0,1]$, $n$ number of pixels in the image and $\sigma_{\veta}$ is the noise level, which is assumed to be known. 
%Finally, $\ell$ is chosen as the least square functional $\ell(\vx,\vb) = \frac{1}{2}\lVert \vA \vx - \vb \rVert^{2}_{2}$.
Since $\rho_{\textsf{\tiny RED}}(\cdot)$ is continuous and the constraint set is also bounded, the problem \eqref{eq:constrained} has at least one solution by the Weierstrass theorem. The problem \eqref{eq:constrained} can be equivalently reformulated into the following form
\begin{align} \label{eq:unconstrained}
    &\argmin{\vx,\vt, \vr \in \mathbb{R}^{n}}{\rho_{\textsf{\tiny RED}}(\vt)} + \iota_{B_{\delta}}(\vr),\, \mbox{ s.t } \ \vr = \vA \vx - \vb, \vx=\vt ,
    %\nonumber&\text{subject to}\quad  \vr = \vA \vx - \vb,\quad \vx=\vt ,
\end{align}
where $B_{\delta}:= \lbrace \vr \in \mathbb{R}^{n} \ \lvert \ \lVert \vr \rVert_{2}^2 \leq \delta \rbrace$ and  $\iota_{B_{\delta}}$ denotes the indicator function of the set $B_{\delta}$. The corresponding Augmented Lagrangian function is given by
\begin{equation}\label{eq:augmented_Lagrangian_constrained}
    \begin{split}
    L(&\vx,\vr, \vt;\vlambda_\vt, \vlambda_\vr) =\, \rho_{\textsf{\tiny RED}}(\vt)+i_{B_\delta}(\vr)\\
    &+\frac{\beta_{\vr}}{2}\lVert \vA \vx - \vb - \vr + \frac{\vlambda_\vr}{\beta_\vr} \rVert_{2}^{2}\\
    &+ \frac{\beta_{\vt}}{2}\left\lVert\vx-\vt+\frac{\vlambda_\vt}{\beta_\vt}\right\rVert_2^2 - \frac12\lVert\vlambda_\vr\rVert^2_2 - \frac12\lVert\vlambda_\vt\rVert^2_2,
    \end{split}
\end{equation}
where $\vlambda_\vr$ and $\vlambda_\vt$ are the Lagrange multipliers, and $\beta_\vr$ and $\beta_\vt$ are real positive penalties. 
The optimization problem involving the Augmented Lagrangian can then solved via the Alternating Direction Method of Multipliers (ADMM) \cite{magnusson2015convergence,wang2019global,cascarano2021combining}. Algorithm \ref{al:redpen} summarizes the ADMM method.

\begin{algorithm}[htbp]
\caption{Constrained RED Approach (CRED)}
\label{al:redpen}
\begin{algorithmic}[1]\small
\STATE{Set $\delta$, $\vx^{0}, \vt^{0}=0$ and $\vr^{0}=\vb-\vA\vx^{0}$, select $\beta_\vu,\, \beta_\vr>0$ and initialise $\vlambda_\vr^{0}, \lambda_\vt^{0}$}
\FOR{$k=0,1,\dots$}
\STATE $\vx^{k+1/2} = \ds\frac{\beta_\vr}{\beta_\vt}\vA^\top\lp \vb +\vr^k -\frac{\vlambda_\vr^k}{\beta_\vr}\rp + \lp\vt^k-\frac{\vlambda_\vt^k}{\beta_\vt}\rp$

\medskip

\STATE $\vx^{k+1} \gets \ds\lp\frac{\beta_\vr}{\beta_\vt}\vA^\top\vA+\vI\rp^{-1} \vx^{k+1/2}$

\medskip

\STATE $\vt^{k+1} = \ds\frac{1}{1+\beta_\vt}f(\vt^k) +\frac{\beta_\vt}{1+\beta_\vt}\lp\vx^{k+1}+\frac{\vlambda_\vt^k}{\beta_\vt}\rp$

\medskip

\STATE $\vr^{k+1} = \ds\proj_{\mathcal{B}_\delta}\lp\vA\vx^{k+1}-\vb+\frac{\vlambda_\vr^k}{\beta_\vr}\rp$

\medskip

\STATE $\vlambda_\vr^{k+1} \gets \vlambda_\vr^k +\beta_\vr(\vr^{k+1}  -\vA\vx^{k+1} +\vb)$

\medskip

\STATE $\vlambda_\vt^{k+1} \gets \vlambda_\vt^k +\beta_\vt\lp-\vx^{k+1} +\vt^{k+1}\rp $
\ENDFOR
\end{algorithmic}
\end{algorithm}
\begin{rem}
    It is worth highlighting the following points relative to Algorithm \ref{al:redpen}. In Line 1 the definition of $\delta$ requires the knowledge of $\sigma_{\veta}$. This is not limiting since a good estimation of the noise standard deviation can often be obtained \cite{immerkaer1996fast}. In lines 3-4 the  optimality conditions on the $\ell_{2}-\ell_{2}$ subproblem are used. This subproblem can often be solved using FFT by imposing periodic boundary conditions. In line 5: we adopt the fixed-point strategy, by zeroing the derivative as in \cite{mataev2019deepred}. In line 6: $\proj_{\mathcal{B}_\delta}$ refers to the projection onto $\mathcal{B}_\delta$. In the implementation, we adopt the variable change $\lambda \sim \frac{\lambda}{\beta}$ for simplifying the notation. Moreover, due to the convexity of $\rho_{\textsf{\tiny RED}}$, the algorithm converges to the minimum of \eqref{eq:constrained}.
\end{rem}
\section{Numerical Results}
\label{sec:numexp}

\subsection{Settings, evaluation metrics, and baseline methods}

We focus on the task of image deblurring with AWGN. Therefore, $\vA$ in \eqref{inverse-problem} is a Gaussian blurring operator of standard deviation $\sigma_{\vA}$.

We simulate blurry and noisy data by applying the image formation model \eqref{inverse-problem} to the images from Set5 \cite{bevilacqua2012low} and Set24 \cite{kodim} referred to as ground truths (GTs). We compare our method with two baselines: the original RED formulation \cite{Romano17} solved using ADMM and its more recent variant RED-PRO \cite{Elad21} solved via gradient-descent. We investigate their behaviour with respect to the choice of some hyperparameters: we focus on the role of $\lambda$ and $(\alpha,\mu)$ for RED and RED-PRO, respectively. The former is the regularization parameter, the couple $(\alpha,\mu)$ represent the strength of the denoiser and the starting steplength, respectively.
Concerning CRED, we consider penalties such that $\beta_{\vr}^{k+1}=\gamma\,\beta_{\vr}^{k}$ and $\beta_{\vt}^{k+1}=~\gamma\,\beta_{\vt}^{k}$. 
In the following, we comment about the choice of $\gamma$ and $\delta$.

For all the methods, we consider the relative difference of the iterates with tolerance equals $10^{-4}$ as the stopping criterion. We set the maximum number of iterations to 200.

We consider the set of CNN based denoisers introduced in \cite{zhang2017learning}. In order to evaluate the influence of the strength of the denoiser, we choose the ones trained for removing Gaussian noise of standard deviation equal to 16 and 30, which are referred to as $\vD_{1}$ and $\vD_{2}$, respectively.

To assess the quality of the restored images we consider the PSNR and SSIM metrics. Moreover, we point out that from a theoretical perspective, given a ground truth image $\vx$ and its blurred and noisy simulated data $\vb$, $\sigma_{\vx}:= \frac{\lVert \vA \vx - \vb \lVert_{2}}{\sqrt{n-1}}$ represents the unbiased estimator of $\sigma_{\veta}$. 
For this reason we consider as valuable metric the comparison between the real noise standard deviation 
$\sigma_{\veta}$ and $\sigma_{\vx^{\ast}}$ where $\vx^{\ast}$ refers to the output of the algorithms.

\subsection{On the choice of the threshold}

In this section, we consider Set5, set $\sigma_{\vA}=1$, $\sigma_{\veta}=25$, $\gamma=1.01$, and use only $\vD_{1}$ as the denoiser. 
We inspect the influence of $\delta=\tau \, \sqrt{n} \, \sigma_{\veta}$, with $\tau \in [0,1]$ on the restored images. We consider two different scenarios: in the first we assume $\sigma_{\veta}$ is known (idealized scenario), whereas in the second we assume only an estimate $\overline{\sigma_{\veta}}$ is provided (realistic scenario). 
The latter is computed as in \cite{immerkaer1996fast}. 

\begin{figure}[htbp]
	\centering
	\subfloat[(a) Plot of $\sigma_{x^{\ast}}(\tau)$.  \label{fig:1_a}]{
	\includegraphics[width=.35\textwidth]{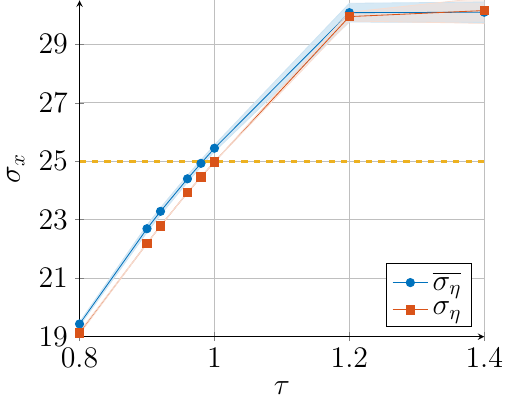}}%\hfill
	\subfloat[(b) Plot of PSNR($\tau$).  \label{fig:1_b}]{
	\includegraphics[width=.35\textwidth]{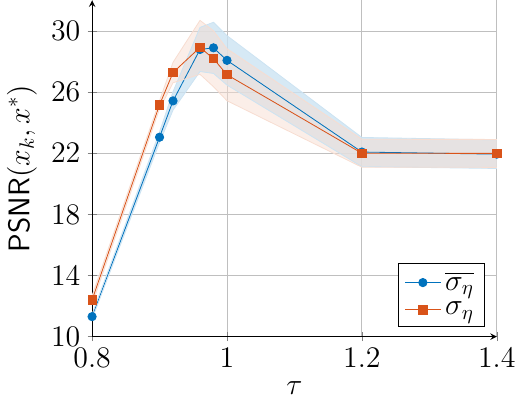}}\hfill
    \caption{Distribution of $\sigma_{\vx^{\ast}}$ (a) and PSNR (b) by varying $\tau$ for idealized and realistic scenarios on the whole Set5.} %over all the images in Set5.} %The continuous lines represent the mean and the shadows represent the standard deviation of the distribution of the  $\sigma_{x^{\ast}}$ and PSNR values over all the images in Set5.}
    \label{fig:1}
\end{figure} 

In Figure \ref{fig:1_a} and Figure~\ref{fig:1_b} the continuous lines show the mean of the distribution of $\sigma_{\vx^{\ast}}$ and PSNR as function of $\tau$. Moreover, we shade the region spanned by the standard deviation of their distributions. The orange and blue lines represent the idealized and realistic scenarios, respectively. The yellow dashed line represents the noise standard deviation $\sigma_{\veta}$ affecting the simulated data.

As expected, in the idealized scenario, we can observe that $\sigma_{\vx^{\ast}}$ approximates $\sigma_{\veta}$ when $\tau=1$ for all images in Set5. Conversely, when only an estimate $\overline{\sigma_{\veta}}$ is given we observe that the best approximation of $\sigma_{\veta}$ is reached when $\tau=0.98$. We hypothesize that this is due to the fact that the algorithm in \cite{immerkaer1996fast} tends to overestimate the noise level in the simulated data $\vb$. Moreover, Figure \ref{fig:1_a} shows that small values of $\tau$ underestimate, whereas large values of $\tau$ overestimate $\sigma_{\veta}$. In Figure \ref{fig:1_b} we show the behaviour of the PSNR metric by varying $\tau$. We obtain comparable performances in terms of PSNR for both scenarios.

For all the following experiments we set $\tau=0.98$ and we estimate $\sigma_{\veta}$ by \cite{immerkaer1996fast}. 
In Table \ref{tab:1} we consider different level of degradations $(\sigma_{\vA}=1, \sigma_{\veta}= 15,25,35,50)$ and we report the mean of the relative errors (RE) between $\sigma_{\vx^{\ast}}$ and $\sigma_{\veta}$ for all the images in Set5.
We observe that the mean of the relative errors is small (less than $0.5 \%$) while changing the level of degradation.

\begin{table}
\centering
\caption{Mean relative errors between $\sigma_{\vx^{\ast}}$ and $\sigma_{\veta}$ for all the images in Set5.}
\begin{tabular}{c|c|c|c|c}
Metric  & $\sigma_{\veta} = 15$ & $\sigma_{\veta} = 25$ & $\sigma_{\veta} = 35$ & $\sigma_{\veta} = 50$ \\
\toprule
\text{RE}($\sigma_{\veta},\sigma_{\vx^{\ast}}$) & 0.0050 & 0.0045 & 0.0044 & 0.0043\\
\end{tabular}
\label{tab:1}
\end{table}

\subsection{On the choice of the ADMM penalties and denoiser \label{subsec:penalties}}

In this section we consider the sole \textit{Butterfly} image by setting $\sigma_{\vA}=1.2$ and $\sigma_{\veta}=30$. We investigate the stability of CRED with respect to the choice of $\gamma$, $\beta_{\vr}^{0}$ and $\beta_{\vt}^{0}$. We point out that in the experiments we did not observe any significant difference when choosing different images.

In Table \ref{tab:gamma} we report the mean values of the distribution of PSNR and SSIM obtained by choosing $\lp\beta_{\vr}^{0}, \beta_{\vt}^{0}\rp \in \{0.2, 0.4, 0.6, 0.8, 1\}^2$ for different $\gamma$. 
In this experiment we consider $\vD_{1}$ as denoiser. 
As a general comment, we observe that we can reach similar performances in terms of PSNR and SSIM for the considered values of $\gamma$. %For the following experiments we set $\gamma=1.01$.

Figure \ref{fig:2_c} depicts the distribution of the PSNR values when setting $\gamma=1.01$ for the CNN denoisers $\vD_{1}$ (orange) and $\vD_{2}$ (blue). 
CRED appears stable regardless the considered choices of $\beta_{\vr}^{0}$ and $\beta_{\vt}^{0}$, and moreover, it seems that the selection of CNN denoisers has a minimal impact on the overall performances.
In the following experiments we fix $\beta_{\vr}^{0}=1$, $\beta_{\vt}^{0}=1$ and $\gamma=1.01$. 

\begin{table}
    \centering
    \caption{Mean values of PSNR and SSIM distributions over different configurations of $\lp\beta_{\vr}^{0}, \beta_{\vt}^{0}\rp$ for different $\gamma$.}
\begin{tabular}{|c|c|c|c|c|}
Metric & $\gamma=1$ & $\gamma=1.01$ & $\gamma=1.05$ \\
\toprule
PSNR & 24.6605 & 24.7312 & 24.4178 \\\midrule
SSIM &  0.9342 & 0.9324  & 0.9219  \\
\end{tabular}

    \label{tab:gamma}
\end{table}

\subsection{Comparisons with RED and RED-PRO}

\begin{figure}[t]%[htpb]
\newcommand\factor{0.30}
	\centering
 \subfloat[(a) RED  \label{fig:2_a}]{\includegraphics[width=\factor\textwidth]{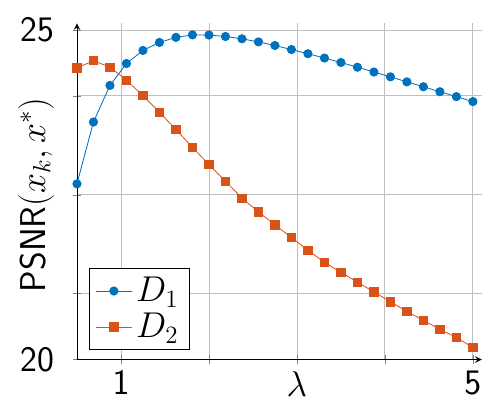}}\hfill
\subfloat[(b) REDPRO  \label{fig:2_b}]{\includegraphics[width=\factor\textwidth]{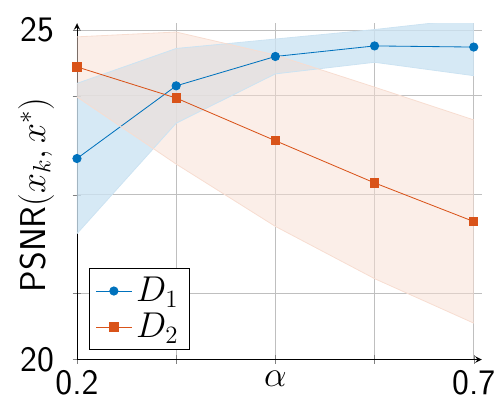}}
\subfloat[(c) CRED  \label{fig:2_c}]{\includegraphics[width=\factor\textwidth]{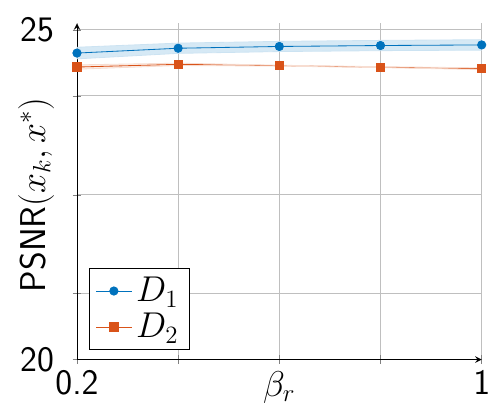}}\hfill 
    \caption{PSNR behaviour by varying  $\lambda$ for RED (a), $(\alpha,\mu)$ and $\mu$ for RED-PRO (b), $(\beta_{\vr},\beta_{\vt})$ for CRED (c). In (b) for each $\alpha$ we present the PSNR distribution (average $\pm$ std) wrt $\mu$; in (c) for each $\beta_\vr$ we present the PSNR distribution (average $\pm$ std) wrt $\beta_\vt$.}
    \label{fig:2}
\end{figure}    

In this section, we compare CRED, RED and RED-PRO in terms of stability with respect to their hyperparameters and reconstruction metric on the Set5 and Set24.

In order to compare the stability of RED and RED-PRO with respect to CRED, we consider the simulated \textit{Butterfly} image used in the Subsection \ref{subsec:penalties}.  
For the original RED algorithm, we sample 25 different values of $\lambda \in [0.5, 5]$. 
Concerning RED-PRO, we consider these 25 different configurations $\lp\alpha, \mu\rp  \in \{0.3, 0.4, 0.5, 0.6, 0.7\}^2$. 

In Figure \ref{fig:2_a} we plot the PSNR behaviour of RED by varying $\lambda$, whereas in Figure \ref{fig:2_b} we report the PSNR distribution of RED-PRO for different configurations $\lp\alpha, \mu\rp$. 
The blue lines represent the case where the denoiser $\vD_{1}$ is used, whereas the orange lines represent the case where the denoiser $\vD_{2}$ is used. 

By comparing these results with the ones reported in Figure \ref{fig:2_c}, it is evident how our CRED looks more stable with respect to the choice of his hyperparameters.
Finally, we observe that for RED and RED-PRO the configuration of parameters maximising the PSNR changes when considering different denoisers. We stress that the same conclusions apply when considering different images.

In order to compare the reconstruction metric performances we consider all the images from Set5 and Set24 for different degradation levels.
In Table \ref{tab:3} we report the mean PSNR and SSIM values. For the competing methods RED and RED-PRO the hyperparameters have been estimated in order to minimize the difference between $\sigma_{\veta}$ and $\sigma_{\vx^{\ast}}$. 
We notice that in terms of the considered metrics CRED performs as well as RED and RED-PRO. However, we remark that it does not require any parameter tuning.   

\begin{table}[hbtp]
\centering
\caption{Mean values of PSNR and SSIM for the images in Set5 and Set24 by varying the degradation levels.}
\begin{tabular}{l|w{c}{0.05\textwidth}w{c}{0.04\textwidth}w{c}{0.04\textwidth}|w{c}{0.05\textwidth}w{c}{0.04\textwidth}w{c}{0.04\textwidth}}
%\hline
& \multicolumn{3}{c|}{Set5 ($\sigma_{\vA}=1, \sigma_{\eta}=15$)} & \multicolumn{3}{c}{Set24 ($\sigma_{\vA}=1.2, \sigma_{\eta}=25$)}\\
\toprule
Metric & RED-PRO & RED & CRED & RED-PRO & RED & CRED \\
\toprule
PSNR & 30.04 & 29.61 & 30.05 & 26.70 & 26.29 & 26.95 \\
\midrule
SSIM & 0.91 & 0.90 & 0.91 & 0.77 & 0.76 & 0.77 \\
\end{tabular}
\label{tab:3}
\end{table}

\begin{figure}
	\centering
 \subfloat[GT \label{fig:GT_butterfly}]{
	\scalebox{0.27}{
	\begin{tikzpicture}
	\begin{scope}[spy using outlines={rectangle,blue,magnification=5,size=5cm}]
	\node [name=c]{	\includegraphics[height=4cm]{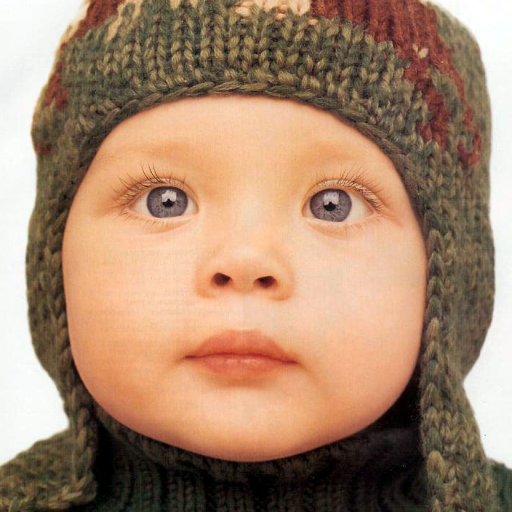}};
        \spy on (-0.1,-0.8) in node [name=c]  at (0,-2.51);
        \spy on (0.7,0.5) in node [name=c]  at (0,2.51);
        %\node [name=c1]{\includegraphics[height=4cm]{images/fig3/Baby_orig.png}};
	\end{scope}
	\end{tikzpicture}}
	}
 \subfloat[Degraded \label{fig:degraded_butterfly}]{
	\scalebox{0.27}{
	\begin{tikzpicture}
	\begin{scope}[spy using outlines={rectangle,blue,magnification=5,size=5cm}]
	\node [name=c]{	\includegraphics[height=4cm]{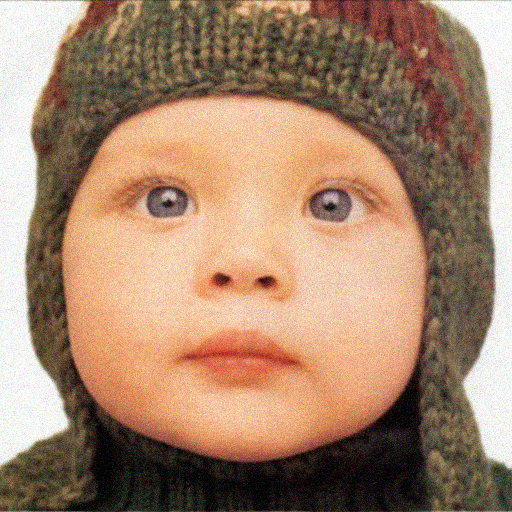}};
	 \spy on (-0.1,-0.8) in node [name=c]  at (0,-2.51);
        \spy on (0.7,0.5) in node [name=c]  at (0,2.51);
	\end{scope}
	\end{tikzpicture}}
	}\subfloat[RED \label{fig:RED_butterfly}]{
	\scalebox{0.27}{
	\begin{tikzpicture}
	\begin{scope}[spy using outlines={rectangle,blue,magnification=5,size=5cm}]
	\node [name=c]{	\includegraphics[height=4cm]{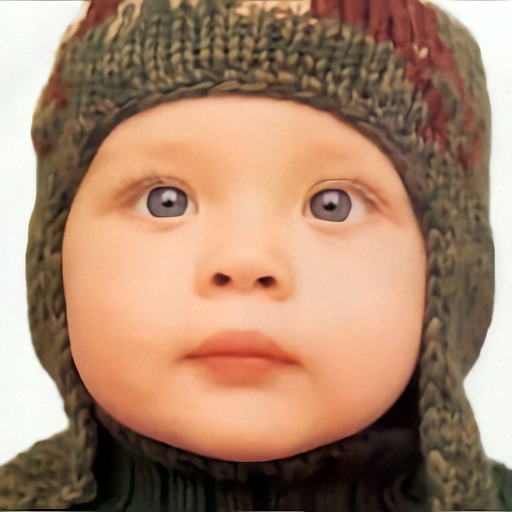}};
	 \spy on (-0.1,-0.8) in node [name=c]  at (0,-2.51);
        \spy on (0.7,0.5) in node [name=c]  at (0,2.51);
	\end{scope}
	\end{tikzpicture}}
	}\subfloat[RED-PRO\label{fig:RED-PRO_butterfly}]{
	\scalebox{0.27}{
	\begin{tikzpicture}
	\begin{scope}[spy using outlines={rectangle,blue,magnification=5,size=5cm}]
	\node [name=c]{	\includegraphics[height=4cm]{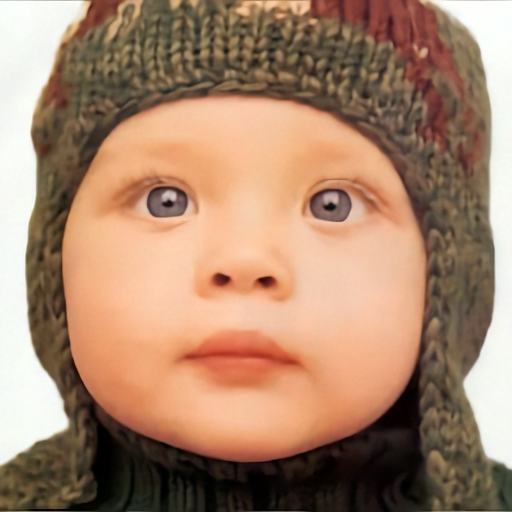}};
	 \spy on (-0.1,-0.8) in node [name=c]  at (0,-2.51);
        \spy on (0.7,0.5) in node [name=c]  at (0,2.51);
	\end{scope}
	\end{tikzpicture}}
	}\subfloat[CRED \label{fig:CRED_butterfly}]{
	\scalebox{0.27}{
	\begin{tikzpicture}
	\begin{scope}[spy using outlines={rectangle,blue,magnification=5,size=5cm}]
	\node [name=c]{	\includegraphics[height=4cm]{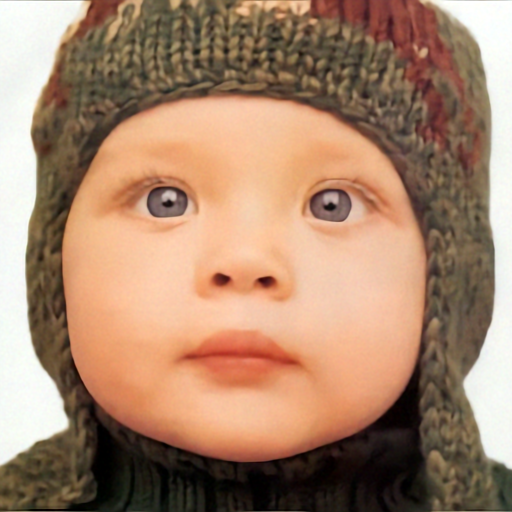}};
	 \spy on (-0.1,-0.8) in node [name=c]  at (0,-2.51);
        \spy on (0.7,0.5) in node [name=c]  at (0,2.51);
	\end{scope}
	\end{tikzpicture}}
	}
	\caption{Restoration of the \emph{Baby} image. From left to right: two close-ups of ground truth, degraded image, RED, RED-PRO, and CRED.}
	\label{fig:3}
    \end{figure}

In Figure \ref{fig:GT_butterfly} and Figure \ref{fig:degraded_butterfly} we report the GT and the degraded close-ups of \textit{Baby} from Set5. Moreover, in Figures \ref{fig:RED_butterfly}, \ref{fig:RED-PRO_butterfly} and \ref{fig:CRED_butterfly} we depict the restored close-ups by RED, RED-PRO and CRED. 
In terms of visual quality there are no relevant differences between the restored images however CRED seems to slightly reproduce more clearly image details than RED and RED-PRO.

     %We present zoomed images for better visualization.  Note the competitive visual performance of CRED.

\section{Conclusions} \label{sec:concl}

This letter presents a novel constrained formulation of the popular RED method that forces the minimum of the regularization functional to satisfy a discrepancy-based bound for a given threshold serving as regularization parameter.  Our formulation is then solved within the ADMM framework and the overall approach results in a simple yet effective method for image restoration, which is called CRED through the letter. Defining the threshold requires an estimate of the standard deviation of the noise affecting the data which can be provided as described in \cite{immerkaer1996fast} thus eliminating the need for extensive parameter estimation. The key point of CRED is its superior stability and robustness with respect to both model and algorithm hyperparameters which is assessed through several comparisons with the original RED and its variant RED-PRO methods. Furthermore, the experiments conducted show that in terms of PSNR and SSIM metrics CRED performs as well as, if not better, when compared to both original RED and RED-PRO. Finally, performances, stability and robustness make CRED a promising choice for image restoration.

\bibliographystyle{IEEEtran}
\bibliography{biblio.bib}
\end{document}